\documentclass[11pt]{amsart}

\usepackage{amsmath}
\usepackage{amssymb}

\newtheorem{theorem}{Theorem}[section]
\newtheorem{claim}[theorem]{Claim}
\newtheorem{mclaim}[theorem]{Main Claim}

\newtheorem{corollary}[theorem]{Corollary}

\theoremstyle{definition}
\newtheorem{definition}[theorem]{Definition}

\newtheorem{question}[theorem]{Question}

\theoremstyle{remark}
\newtheorem{remark}[theorem]{Remark}

\newcount\skewfactor
\def\mathunderaccent#1#2 {\let\theaccent#1\skewfactor#2
\mathpalette\putaccentunder}
\def\putaccentunder#1#2{\oalign{$#1#2$\crcr\hidewidth
\vbox to.2ex{\hbox{$#1\skew\skewfactor\theaccent{}$}\vss}\hidewidth}}
\def\name{\mathunderaccent\tilde-3 }

\def\smallbox#1{\leavevmode\thinspace\hbox{\vrule\vtop{\vbox
   {\hrule\kern1pt\hbox{\vphantom{\tt/}\thinspace{\tt#1}\thinspace}}
   \kern1pt\hrule}\vrule}\thinspace}

\newcommand{\cf}{{\rm cf}}

\def\qedref#1{$\qed_{\reforiginal{#1}}$}

\title{Random reals and polarized colorings}
\author{Shimon Garti}
\address{Institute of Mathematics,
 The Hebrew University of Jerusalem,
 Jerusalem 91904, Israel}
\email{shimon.garty@mail.huji.ac.il}

\author{Saharon Shelah}
\address{Institute of Mathematics
 The Hebrew University of Jerusalem
 Jerusalem 91904, Israel
 and  Department of Mathematics
 Rutgers University
 New Brunswick, NJ 08854, USA}
\email{shelah@math.huji.ac.il}
\urladdr{http://www.math.rutgers.edu/\char`\~shelah}

\thanks{The second author is supported by the European Research Council, Grant 338821. This is publication 1047 of the second author}

\begin{document}
\let\labeloriginal\label
\let\reforiginal\ref

\begin{abstract}
We analyze the strong polarized partition relation with respect to several cardinal characteristics and forcing notions of the reals. We prove that random reals (as well as the existence of real-valued measurable cardinals) yield downward negative polarized relations.
\end{abstract}

\maketitle

\newpage

\section{Introduction}

This paper focuses on two cardinal characteristics of the continuum, the reaping number $\mathfrak{r}$ and the splitting number $\mathfrak{s}$. Let us commence with the basic definitions of these invariants:

\begin{definition}
\label{rrrr} The reaping number.
\begin{enumerate}
\item [$(\aleph)$] Suppose $B\in[\omega]^\omega$ and $S\subseteq\omega$. $S$ splits $B$ if $|S\cap B|=|(\omega\setminus S)\cap B|=\aleph_0$.
\item [$(\beth)$] $\{T_\alpha:\alpha<\kappa\}$ is an unreaped family if there is no single $S\in[\omega]^\omega$ so that $S$ splits $T_\alpha$ for every $\alpha<\kappa$.
\item [$(\gimel)$] The reaping number $\mathfrak{r}$ is the minimal cardinality of an unreaped family.
\item [$(\daleth)$] $\mathfrak{r}_\sigma$ is the minimal cardinality of a collection $\mathcal{R}\subseteq[\omega]^\omega$ which is not splitted by $\omega$-many sets.
\end{enumerate}
\end{definition}

The dual of the reaping number is the splitting number.
Recall:

\begin{definition}
\label{ssss} The splitting number.
\begin{enumerate}
\item [$(\aleph)$] Suppose $B\in[\omega]^\omega$ and $S\subseteq\omega$. $S$ splits $B$ if $|S\cap B|=|(\omega\setminus S)\cap B|=\aleph_0$.
\item [$(\beth)$] $\{S_\alpha:\alpha<\kappa\}$ is a splitting family in $\omega$ if for every $B\in[\omega]^\omega$ there exists an ordinal $\alpha<\kappa$ so that $S_\alpha$ splits $B$.
\item [$(\gimel)$] The splitting number $\mathfrak{s}$ is the minimal cardinality of a splitting family in $\omega$.
\end{enumerate}
\end{definition}

In the present paper, combinatorial arguments serve in the area of cardinal invariants of the continuum. The main tool is the following concept.
If $\lambda\geq\kappa$ are infinite cardinals then the strong polarized relation $\binom{\lambda}{\kappa}\rightarrow\binom{\lambda}{\kappa}^{1,1}_2$ means that for every $c:\lambda\times\kappa\rightarrow 2$ there are $A\in[\lambda]^\lambda, B\in[\kappa]^\kappa$ such that $c\upharpoonright(A\times B)$ is constant.
We shall make use of the following theorem from \cite{MR3000439}:

\begin{claim}
\label{belowssss} Strong polarized relations below the splitting number. \newline 
Assume $\kappa<\mathfrak{s}$.\newline 
The positive relation $\binom{\kappa}{\omega}\rightarrow \binom{\kappa}{\omega}^{1,1}_2$ holds iff $\cf(\kappa)>\aleph_0$.
\end{claim}

\hfill \qedref{belowssss}

In a way, the splitting number $\mathfrak{s}$ is a natural point for proving downward positive relations. The dual notion of the reaping number $\mathfrak{r}$ is a natural point for upward positive relations, as shown in the following theorem from \cite{MR3201820}:

\begin{theorem}
\label{aboverrrr} Strong polarized relations above the reaping number. \newline 
Assume $\mathfrak{r}<\kappa\leq\mathfrak{c}$. \newline 
If $\mathfrak{r}<\cf(\kappa)$ then $\binom{\kappa}{\omega}\rightarrow \binom{\kappa}{\omega}^{1,1}_2$.
\end{theorem} 

\hfill \qedref{aboverrrr}

Observe that the requirement about the cofinality of $\kappa$ is stronger in this theorem, and we do not have full knowledge when $\aleph_0<\cf(\kappa)\leq\mathfrak{r}$, see below. However, some kind of duality is reflected in this theorem. The keypoint for the results of this paper is that we can prove also a \emph{negative} downward theorem with respect to $\mathfrak{r}$. It enables us to supply a positive answer to the following open problem from \cite{MR3509813} (Problem 3.19 there):

\begin{question}
\label{from1012} Suppose $\aleph_1<\kappa=\cf(\kappa)<\lambda=\mathfrak{c}$ and $\binom{\kappa}{\omega}\rightarrow \binom{\kappa}{\omega}^{1,1}_2$. 
Let $\mathbb{P}$ be L\'evy$(\kappa,\lambda)$. Is it possible that $\binom{\kappa}{\omega}\nrightarrow \binom{\kappa}{\omega}^{1,1}_2$ in ${\rm\bf{V}}^{\mathbb{P}}$?
\end{question}

As we shall see, a wide range of positive answers can be given to the above problem, i.e., for every regular cardinal $\kappa$ we can build a model of ZFC in which the L\'evy collapse destroys the strong polarized relation for $\kappa$.

Another problem from \cite{MR3509813} focuses on the random real forcing. One of the salient properties of random reals is that they are dominated by old reals from the ground model, hence the dominating number $\mathfrak{d}$ is unchanged. Problem 3.9 (in \cite{MR3509813}) asks if one can ruin the positive relation for $\mathfrak{d}$ while iterating random real forcing notions:

\begin{question}
\label{dddd} Assume $\binom{\mathfrak{d}}{\omega}\rightarrow \binom{\mathfrak{d}}{\omega}^{1,1}_2$, and one adds $\lambda$-many random reals (for some $\lambda>\mathfrak{d}$). Does the positive relation $\binom{\mathfrak{d}}{\omega}\rightarrow \binom{\mathfrak{d}}{\omega}^{1,1}_2$ still hold?
\end{question}

As in the former question, we will be able to show that one can force a negative relation for $\mathfrak{d}$ after adding random reals to the universe. This follows, again, from the downward negative relation for $\mathfrak{r}$. Moreover, the implication of negative results upon adding many random reals is wider. This gives rise to the third problem (number 3.12) that we quote:

\begin{question}
\label{rrvm} Assume $\kappa$ is a real-valued measurable cardinal. \newline 
Is it possible that $\binom{\kappa}{\omega}\rightarrow \binom{\kappa}{\omega}^{1,1}_2$?
\end{question}

Concerning this question, we will be able to supply only a partial answer, by proving many negative relations \emph{below} the real-valued measruable cardinal. We comment that these relations would give the interesting corollary that $\mathfrak{s}=\aleph_1$ whenever a real-valued measurable cardinal exists. Likewise, under the additional assumption that $\kappa=\mathfrak{c}$ is real-valued measurable we shall prove that $\binom{\kappa}{\omega}\nrightarrow \binom{\kappa}{\omega}^{1,1}_2$ as well.

Our notation is standard. We shall follow \cite{MR2768685} with respect to cardinal invariants. We employ the Jerusalem forcing notation, so $p\leq q$ means that $q$ is stronger than $p$. In cases where ambiguity lurks around the corner we shall state the pertinent conventions explicitly.
We thank the referee for the careful reading of the manuscript.

\newpage 

\section{Continuum reaping}

Investigating the interplay between cardinal invariants and strong polarized relations, one may ask whether there exists a natural cardinal invariant $\mathfrak{y}$ so that $\binom{\mathfrak{y}}{\omega}\nrightarrow\binom{\mathfrak{y}}{\omega}^{1,1}_2$ always holds. The answer is negative. It is shown in \cite{MR3201820} that for every cardinal invariant, the above polarized relation is independent. Nonetheless, there is a natural characteristic for which the negative relation follows, under the extra assumption of being the continuum.
Our basic result says that if the reaping number and the continuum coincide, then a negative polarized parition relation can be proved dwon to the cofinality of the continuum:

\begin{mclaim}
\label{mc} Negative downward relations. \newline 
If $\mathfrak{r}=\mathfrak{c}$
Then $\binom{\mathfrak{r}}{\omega}\nrightarrow\binom{\mathfrak{r}}{\omega}^{1,1}_2$. \newline  
Moreover, $\binom{\kappa}{\omega}\nrightarrow \binom{\kappa}{\omega}^{1,1}_2$ for every $\kappa\in[\cf(\mathfrak{c}),\mathfrak{c}]$.
\end{mclaim}

\par\noindent\emph{Proof}. \newline 
Enumerate the members of $[\omega]^\omega$ by $\{B_\gamma:\gamma<\mathfrak{r}\}$.
For every $\alpha<\mathfrak{r}$ let $\mathcal{B}_\alpha=\{B_\gamma:\gamma<\alpha\}$. The size of $\mathcal{B}_\alpha$ is less than $\mathfrak{r}$, so we can choose a set $S_\alpha$ which splits all the members of $\mathcal{B}_\alpha$. This is rendered for every $\alpha<\mathfrak{r}$, and gives rise to a coloring $d:\mathfrak{r}\times\omega\rightarrow\{0,1\}$ as follows:

$$
d(\alpha,n)=0 \Leftrightarrow n\in S_\alpha.
$$

We claim that $d$ exemplifies the negative relation $\binom{\mathfrak{r}}{\omega}\nrightarrow\binom{\mathfrak{r}}{\omega}^{1,1}_2$. Indeed, assume $H\in[\mathfrak{r}]^{\mathfrak{r}}$ and $B\in[\omega]^\omega$. Assume toward contradiction that $d\upharpoonright(H\times B)$ is constant. If the constant value is $0$ then $B\subseteq S_\alpha$ for every $\alpha\in H$, and if the constant value is $1$ then $B\subseteq (\omega\setminus S_\alpha)$ for every $\alpha\in H$. 

In any case,
The set $B$ appears in the above enumeration, so $B\equiv B_\gamma$ for some $\gamma<\mathfrak{r}$. Choose an ordinal $\alpha\in H$ so that $\gamma<\alpha$, and notice that $S_\alpha$ splits $B$, a contradiction.

Moreover, if $\kappa\geq\cf(\mathfrak{c})$ then we enumerate the members of $[\omega]^\omega$ as $\{B_\gamma:\gamma<\mathfrak{c}\}$. We choose an increasing and unbounded sequence of ordinals $\langle\alpha_\varepsilon:\varepsilon<\kappa\rangle$ in $\mathfrak{c}$ (repetitions are welcome), and define $\mathcal{B}_\varepsilon=\{B_\gamma:\gamma<\alpha_\varepsilon\}$ for every $\varepsilon<\kappa$. Again, $|\mathcal{B}_\varepsilon|<\mathfrak{r}$ for every $\varepsilon<\kappa$ as $\mathfrak{r}=\mathfrak{c}$. For every $\varepsilon<\kappa$ we choose some $S_\varepsilon\in[\omega]^\omega$ which splits all the members of $\mathcal{B}_\varepsilon$.

The coloring $d:\kappa\times\omega\rightarrow\{0,1\}$ is defined in the same way, i.e., $d(\varepsilon,n)=0 \Leftrightarrow n\in S_\varepsilon$. The same argument shows that $d$ has no monochromatic product of size $\kappa\times\omega$, so we are done.

\hfill \qedref{mc}

\begin{remark}
\label{highrr} For any uncountable cardinal $\theta$ define $\mathfrak{r}_\theta$ as the minimal cardinality of a subset of $[\theta]^\theta$ such that no single $B\in[\theta]^\theta$ splits all the members of this family. One can verify that $\mathfrak{r}_\theta>\theta$ for every infinite cardinal $\theta$, and the main claim holds for every $\theta$ (under the parallel generalized assumption that $\mathfrak{r}_\theta=2^\theta$).
\end{remark}

\hfill \qedref{highrr}

The main claim is optimal in the sense that the assumption $\mathfrak{r}=\mathfrak{c}$ cannot induce a stronger negative relation in ZFC. The closest attempt would be refuting the positive unbalanced relation $\binom{\mathfrak{r}}{\omega}\rightarrow \binom{\mathfrak{r}\ \alpha}{\omega\ \omega}^{1,1}_2$ for every $\alpha<\mathfrak{r}$, but the following claim proves its independence:

\begin{claim}
\label{uunbalanced} Unbalanced negative relations. \newline 
The assumption $\mathfrak{r}=\mathfrak{c}$ is consistent with both $\binom{\mathfrak{r}}{\omega}\rightarrow \binom{\mathfrak{r}\ \alpha}{\omega\ \omega}^{1,1}_2$ for every $\alpha<\mathfrak{r}$ and $\binom{\mathfrak{r}}{\omega}\nrightarrow \binom{\mathfrak{r}\ \alpha}{\omega\ \omega}^{1,1}_2$ for some $\alpha<\mathfrak{r}$.
\end{claim}

\par\noindent\emph{Proof}. \newline 
For the positive direction force $\mathfrak{p}=\mathfrak{c}$, in which case $\mathfrak{r}=\mathfrak{c}$ as well. However, the relation $\binom{\mathfrak{p}}{\omega}\rightarrow \binom{\mathfrak{p}\ \alpha}{\omega\ \omega}^{1,1}_2$ for every $\alpha<\mathfrak{p}$ is established in \cite{MR2444279} and holds in ZFC, so $\binom{\mathfrak{r}}{\omega}\rightarrow \binom{\mathfrak{r}\ \alpha}{\omega\ \omega}^{1,1}_2$ for every $\alpha<\mathfrak{r}$. Observe that $\mathfrak{r}$ is a regular cardinal in such models.

For the negative direction choose any $\lambda=\lambda^{\aleph_0}$ such that $\lambda>\aleph_1$. Let $\mathbb{Q}$ be a finite support iteration of adding $\lambda$-many Cohen reals. It is known that ${\rm \bf V}^{\mathbb{Q}}\models\binom{\mu}{\omega}\nrightarrow \binom{\omega_1}{\omega}^{1,1}_2$ for every $\mu\in(\aleph_0,\lambda]$, as shown in \cite{MR3000439}, Remark 2.4. In particular, it holds for $\mu=\mathfrak{c}$. As $\mathfrak{r}=\mathfrak{c}$ in this generic extension and $\lambda>\aleph_1$ we have the consistency of the negative direction.

\hfill \qedref{uunbalanced}

We turn back to the balanced relation. In the case of a regular continuum, we can characterize now the strong polarized relation for $\mathfrak{c}$ as follows:

\begin{corollary}
\label{ccc} Assume $\mathfrak{c}$ is a regular cardinal. \newline 
Then $\binom{\mathfrak{c}}{\omega}\rightarrow\binom{\mathfrak{c}}{\omega}^{1,1}_2$ iff $\mathfrak{r}<\mathfrak{c}$.
\end{corollary}

\par\noindent\emph{Proof}. \newline 
If $\mathfrak{r}<\mathfrak{c}$ then Theorem \ref{aboverrrr} gives the positive direction of $\binom{\mathfrak{c}}{\omega}\rightarrow\binom{\mathfrak{c}}{\omega}^{1,1}_2$, since $\mathfrak{c}$ is a regular cardinal. If $\mathfrak{r}=\mathfrak{c}$ then Claim \ref{mc} gives the negative relation, so the proof is accomplished.

\hfill \qedref{ccc}

We employ the above corollary in the proof of the following theorem:

\begin{theorem}
\label{mt} L\'evy collapse and random reals. \newline 
Suppose $\kappa$ is an uncountable regular cardinal. \newline 
For every $\lambda=\cf(\lambda)>\kappa$ there is a model of ZFC in which $\mathfrak{c}=\lambda, \binom{\kappa}{\omega}\rightarrow \binom{\kappa}{\omega}^{1,1}_2$ and if $\mathbb{P}=$L\'evy$(\kappa,\lambda)$ then ${\rm\bf{V}}^{\mathbb{P}}\models \binom{\kappa}{\omega}\nrightarrow \binom{\kappa}{\omega}^{1,1}_2$. \newline 
Likewise, it is consistent that $\binom{\mathfrak{d}}{\omega}\rightarrow \binom{\mathfrak{d}}{\omega}^{1,1}_2$ in the ground model, and after adding $\lambda$-many random reals for some $\lambda>\mathfrak{d}$ we have $\binom{\mathfrak{d}}{\omega}\nrightarrow \binom{\mathfrak{d}}{\omega}^{1,1}_2$.
\end{theorem}

\par\noindent\emph{Proof}. \newline 
If $\kappa=\aleph_1$ then the theorem follows from the fact that $2^{\aleph_0}=\aleph_1$ implies $\binom{\aleph_1}{\aleph_0}\nrightarrow \binom{\aleph_1}{\aleph_0}^{1,1}_2$ (as proved in \cite{MR0202613}). So assume that $\kappa>\aleph_1$. 
We begin with MA + $2^{\aleph_0}=\lambda$. In this case, $\mathfrak{s}=\lambda$ as well, so $\binom{\kappa}{\omega}\rightarrow \binom{\kappa}{\omega}^{1,1}_2$ by Claim \ref{belowssss}. Observe also that $\mathfrak{r}=\lambda$. We claim that after forcing with $\mathbb{P}$ we will get the negative relation $\binom{\kappa}{\omega}\nrightarrow \binom{\kappa}{\omega}^{1,1}_2$.

Indeed, $\mathfrak{r}=\kappa$ in the generic extension. This fact follows from the completeness of $\mathbb{P}$ which ensures that no new sequence of sets of length below $\kappa$ is introduced. The length $\lambda$ of $\mathfrak{r}$-sequences in the old universe is collapsed to $\kappa$, but no $\mathfrak{r}$-family of size less than $\kappa$ appears. Hence $\mathfrak{r}=\mathfrak{c}=\kappa$ in ${\rm\bf{V}}^{\mathbb{P}}$. From Claim \ref{mc} we infer that $\binom{\kappa}{\omega}\nrightarrow \binom{\kappa}{\omega}^{1,1}_2$ as required.

We indicate that if $\mathfrak{r}<\kappa$ in the ground model then the positive relation $\binom{\kappa}{\omega}\rightarrow \binom{\kappa}{\omega}^{1,1}_2$ holds both in the old universe and after the collapse (see Corollary \ref{ccc}), so the opposite situation is also consistent for every regular cardinal $\kappa$ above $\aleph_1$.

For the second assertion, begin with a model in which the positive relation $\binom{\mathfrak{d}}{\omega}\rightarrow \binom{\mathfrak{d}}{\omega}^{1,1}_2$ holds, and $\mathfrak{d}>\aleph_1$. This can be done due to \cite{MR3201820}, based on the model of \cite{MR879489}, upon noticing that $\mathfrak{r}<\mathfrak{d}$ gives the desired result when $\mathfrak{d}$ is a regular cardinal (see Theorem \ref{aboverrrr}).

We choose a large enough singular cardinal $\lambda$ so that $\lambda>\mathfrak{d}$ but $\cf(\lambda)\leq\mathfrak{d}$. By adding $\lambda$-many random reals we blow up the continuum to $\lambda$ but $\mathfrak{d}$ remains in its place. Moreover, $\mathfrak{r}=\lambda$ as well (see, e.g., \cite{MR2768685}). Since $\cf(\lambda)\leq\mathfrak{d}$ we conclude that the negative relation $\binom{\mathfrak{d}}{\omega}\nrightarrow \binom{\mathfrak{d}}{\omega}^{1,1}_2$ holds in the generic extension, so the proof is accomplished.

\hfill \qedref{mt}

Can we incorporate singular cardinals in Corollary \ref{ccc}? A good understanding of polarized relations for singular cardinals above $\mathfrak{r}$ is needed. It is consistent that $\kappa>\cf(\kappa)=\mathfrak{r}$ and $\binom{\kappa}{\omega}\rightarrow \binom{\kappa}{\omega}^{1,1}_2$. The opposite direction is not so clear:

\begin{question}
\label{q1} Is it consistent that $\kappa>\mathfrak{r}, \cf(\kappa)>\aleph_0$ and $\binom{\kappa}{\omega}\nrightarrow \binom{\kappa}{\omega}^{1,1}_2$?
\end{question}

The negative downward theorem below $\mathfrak{r}$ (under the assumption that $\mathfrak{r}=\mathfrak{c}$) can be used also for a surprising relationship between $\mathfrak{r}$ and $\mathfrak{s}$. One of the dividing lines in the realm of cardinal invariants is the distinction between small characteristics (which are bounded by $\cf(\mathfrak{c})$) and large characteristics (which are not bounded by $\cf(\mathfrak{c})$). The distributivity number $\mathfrak{h}$ is a typical example of a small invariant, while $\mathfrak{s}$ is a large invariant. Nevertheless, if $\mathfrak{r}=\mathfrak{c}$ then $\mathfrak{s}$ becomes small:

\begin{theorem}
\label{smallssss} If $\mathfrak{r}=\mathfrak{c}$
then $\mathfrak{s}\leq\cf(\mathfrak{c})$. \newline 
Moreover, if $\mathfrak{r}_\sigma=\mathfrak{c}$ then $\mathfrak{s}\leq\cf(\mathfrak{c})$.
\end{theorem}

\par\noindent\emph{Proof}. \newline 
We prove the first assertion with the aid of the polarized relations, and the second assertion in a direct way. Let $\kappa$ be $\cf(\mathfrak{c})$. By Claim \ref{mc} we have $\binom{\kappa}{\omega}\nrightarrow \binom{\kappa}{\omega}^{1,1}_2$, since $\mathfrak{r}=\mathfrak{c}$. This relation excludes the possibility that $\mathfrak{s}>\cf(\mathfrak{c})$, because in this situation we have $\binom{\kappa}{\omega}\rightarrow \binom{\kappa}{\omega}^{1,1}_2$ since $\kappa$ is an uncountable regular cardinal and due to Claim \ref{belowssss}.

Assume now that $\mathfrak{r}_\sigma=\mathfrak{c}$. Enumerate the members of $[\omega]^\omega$ by $\{B_\gamma:\gamma<\mathfrak{c}\}$, and choose an increasing unbounded sequence of ordinals of the form $\langle\alpha_\varepsilon:\varepsilon<\kappa\rangle$ in $\mathfrak{c}$. For every $\varepsilon<\kappa$ let $\mathcal{B}_\varepsilon$ be $\{B_\gamma:\gamma<\alpha_\varepsilon\}$, and we choose a collection of sets $\{S^\varepsilon_n:n\in\omega\}$ which splits the members of $\mathcal{B}_\varepsilon$.

The collection $\mathcal{F}=\{S^\varepsilon_n:\varepsilon<\kappa, n\in\omega\}$ is a splitting family for $[\omega]^\omega$. 
For this, pick up any $B\in[\omega]^\omega$ and any ordinal $\varepsilon<\kappa$ so that $B\in\mathcal{B}_\varepsilon$. By the choice of $\{S_n^\varepsilon: n\in\omega\}$ there is a set $S_n^\varepsilon$ which splits $B$. But $S_n^\varepsilon\in\mathcal{F}$, hence $\mathcal{F}$ is a splitting family.
Consequently, $\mathfrak{s}\leq|\mathcal{F}|\leq\cf(\mathfrak{c})$, so we are done.

\hfill \qedref{smallssss}

The above results raise some natural problems. We phrase a couple of them:

\begin{question}
\label{qqqq} Small cofinality above $\mathfrak{r}$.
\begin{enumerate}
\item [$(\alpha)$] Assume $\mathfrak{r}<\kappa\leq\mathfrak{c}$ and $\cf(\kappa)\leq\mathfrak{r}$. Is it possible that $\binom{\kappa}{\omega}\nrightarrow\binom{\kappa}{\omega}^{1,1}_2$? In particular, is it possible for $\kappa=\mathfrak{c}$?
\item [$(\beta)$] Assume $\mathfrak{r}_\sigma=\mathfrak{c}$. Is it provable that $\binom{\mathfrak{r}_\sigma}{\omega}\nrightarrow \binom{\mathfrak{r}_\sigma}{\omega}^{1,1}_2$?
\end{enumerate}
\end{question}

\newpage 

\section{Random reals and real-valued measurability}

In this section we try to analyze the polarized relation under the existence of random reals and in the presence of real-valued measurable cardinals. 
For a general background and notational conventions used below, we refer to \cite{MR2250550}.
We commence with a negative downward spectrum which issues from adding random reals.

\begin{theorem}
\label{aaddingrandoms} Random reals and polarized relations. \newline 
Assume $\kappa>\aleph_0$ and $\mathbb{Q}$ is a forcing notion for adding $\kappa$-many random reals. \newline 
Then $\Vdash_{\mathbb{Q}} \binom{\theta}{\omega}\nrightarrow\binom{\theta}{\omega}^{1,1}_2$ for every $\theta\leq\kappa$, and even $\Vdash_{\mathbb{Q}} \binom{\kappa}{\omega}\nrightarrow\binom{\theta}{\omega}^{1,1}_2$.
\end{theorem}

\par\noindent\emph{Proof}. \newline 
Choose a generic subset $G\subseteq\mathbb{Q}$, and fix a cardinal $\theta\in[\aleph_1,\kappa]$. Let $m$ be the product measure over ${}^\kappa 2$. Recall that $p\in\mathbb{Q}$ iff $p\subseteq{}^\kappa 2$, where $p$ is a Borel set of positive measure, supported by some countable set $u=u_p\in[\kappa]^{\leq\aleph_0}$. 
We indicate that a support of a given condition $p$ is not unique (every countable subset of $\kappa$ which contains a support can serve as well), though a minimal support always exists.
For $p,q\in\mathbb{Q}$ we define $p\leq q$ iff $q\subseteq p$.

Let $\mathcal{F}$ be the set $\{f:f$ is a finite (partial) function from $\kappa$ into $2\}$. For every $f\in\mathcal{F}$ let $\mathcal{O}(f)$ be $\{g\in{}^\kappa 2:f\subset g\}$, denoted also by $({}^\kappa 2)^{[f]}$. By the definition of the product measure, $m(\mathcal{O}(f))=\frac{1}{2^{|f|}}$. In particular, each $\mathcal{O}(f)$ belongs to $\mathbb{Q}$.

We define a $\mathbb{Q}$-name $\name{\eta}$ of a function from $\kappa$ into $\{0,1\}$ by $\name{\eta}\supset f \Leftrightarrow \mathcal{O}(f)\in G$. 
If $\alpha\in\kappa, f=\{\langle\alpha,0\rangle\}$ and $g=\{\langle\alpha,1\rangle\}$ then $\mathcal{A}=\{\mathcal{O}(f),\mathcal{O}(g)\}$ forms a maximal antichain, so exactly one member of $\mathcal{A}$ belongs to $G$. Hence $\alpha\in{\rm dom}(\name{\eta})$ for every $\alpha<\kappa$. The fact that $\name{\eta}$ is a function follows from the directness of $G$.
This gives rise to the definition of a name $\name{c}$ for a coloring from $\kappa\times\omega$ into $\{0,1\}$ as follows:
$$
\name{c}(\alpha,n)=\name{\eta}(\omega\alpha+n).
$$
We claim that $\name{c}$ exemplifies the negative relation $\binom{\kappa}{\omega}\nrightarrow\binom{\theta}{\omega}^{1,1}_2$.

Assume this is not the case. Pick up a condition $p\in\mathbb{Q}$, a color $\ell\in\{0,1\}$ and names of sets $\name{A}\in[\kappa]^\theta,\name{B}\in [\omega]^{\aleph_0}$ so that $p\Vdash\name{c}\upharpoonright(\name{A}\times \name{B})=\{\ell\}$.

For every $\alpha<\kappa$ we choose a Borel set $B_\alpha\subseteq{}^\kappa 2$ with support $u_\alpha\in[\kappa]^{\leq\aleph_0}$ such that $B_\alpha\subseteq p$ and $B_\alpha$ decides whether $\check{\alpha}$ belongs to $\name{A}$ in the following sense: if $m(B_\alpha)>0$ then $B_\alpha\Vdash\check{\alpha}\in\name{A}$ and if not then $(p-B_\alpha)\Vdash\check{\alpha}\notin\name{A}$.
Similarly, for every $n\in\omega$ we choose a Borel set $B_n\subseteq{}^\kappa 2$ with support $u_n\in[\kappa]^{\leq\aleph_0}$ such that $B_n\subseteq p$ and the following is satisfied: if $m(B_n)>0$ then $B_n\Vdash\check{n}\in\name{B}$ and if not then $(p-B_n)\Vdash\check{n}\notin\name{B}$.

Let $v\in[\kappa]^{\leq\aleph_0}$ be a support of the condition $p$ so that $u_n\subseteq v$ for every $n\in\omega$. We shall force with the part of $\mathbb{Q}$ above $p$ using conditions with the support $v$.

For any ordinal $\alpha\in\kappa\setminus v$ let $r_\alpha=\{\eta\in{}^\kappa 2:\eta\upharpoonright\{\omega\alpha+n:n\in B\}=\ell\}$. Since $\name{B}$ is a name of an unbounded subset of $\omega$, $m(r_\alpha)=0$. Indeed, for every $n\in\omega$ let $\{b_j:j<n\}$ enumerate the first $n$ members of $B$, and let $f_n=\{\langle\omega\alpha+b_j,\ell\rangle:j<n\}$. By definition, $m(\mathcal{O}(f_n))=\frac{1}{2^n}$, and $m(r_\alpha)\leq m(\mathcal{O}(f_n))$ for every $n\in\omega$, so $m(r_\alpha)=0$.

Since $p\Vdash\name{A}\in[\kappa]^\theta$, there exists an ordinal $\alpha$ and a condition $q\geq p$ so that $q\Vdash\check{\alpha}\in\name{A}$ and $\alpha\notin v$. Notice that $q\Vdash\name{c}(\alpha,n)=\ell$ for every $n\in B$. It follows that $q\subseteq r_\alpha$, so $m(q)=0$, which is impossible since $q\in\mathbb{Q}$.

\hfill \qedref{aaddingrandoms}

The effect of adding random reals is sharpened if we assume the existence of a real-valued measurable cardinal. The classical way to introduce such a cardinal is the random real forcing, as proved by Solovay, but this is not the only way. 
Gitik and Shelah, \cite{MR1035887}, introduced a different way to introduce such cardinals, and some of the properties of Solovay's construction are not shared by all real-valued measurable cardinals.
We shall see, however, that the mere existence of a real-valued measurable cardinal entails strong negative relations.

We say that $\kappa$ is real-valued measurable iff there exists an atomless $\kappa$-additive measure over $\kappa$. The requirement of being atomless implies $\kappa\leq 2^{\aleph_0}$, so $\kappa$ is not strongly inaccessible. It is known, however, that such $\kappa$ is weakly inaccessible.
Before embarking on the impact of real-valued meaurable cardinals we need some preliminaries. Let $(X,\Sigma,m)$ be a measure space. The measure algebra associated with it is the Boolean algebra $B=\Sigma/I$ when $I=\{A\subseteq X: m(A)=0\}$. The following belongs to Maharam:

\begin{theorem}
\label{maharamthm} Maharam's Theorem. \newline 
Suppose $m$ is a homogeneous $\sigma$-additive measure on a $\sigma$-complete Boolean algebra $\mathcal{B}$. \newline 
The measure algebra $(\mathcal{B},m)$ is isomorphic to the measure algebra of ${}^\lambda 2$ for some $\lambda$, with the product measure.
\end{theorem}

This fundamental theorem appears in \cite{MR0006595}. 
Let $\kappa$ be real-valued measurable as witnessed by the measure $m$, and let $\mathcal{I}=\{a\subseteq \kappa:m(a)=0\}$. From Maharam's theorem there is some $\mu$ for which $\mathcal{P}(\kappa)/\mathcal{I}$ is isomorphic to the Boolean algebra ${\rm Borel}({}^\mu 2)/\mathcal{J}$, where $\mathcal{J}$ is the ideal of null sets in ${}^\mu 2$.
It has been proved in \cite{MR1035887}, Section 2, that $\mu>\kappa$.

\begin{theorem}
\label{nnegbelowrvm} Negative relations and real-valued measurable cardinals. \newline 
Let $\kappa$ be a real-valued measurable cardinal. \newline 
Then $\binom{\theta}{\omega}\nrightarrow\binom{\aleph_1}{\omega}^{1,1}_2$ for every $\theta<\kappa$.
\end{theorem}

\par\noindent\emph{Proof}. \newline 
Let $m:\mathcal{P}(\kappa)\rightarrow[0,1]_{\mathbb{R}}$ be a measure which exemplifies the fact that $\kappa$ is real-valued measurable. The collection of sets $\mathcal{I}=\{a\subseteq \kappa:m(a)=0\}$ is a $\kappa$-complete ideal over $\kappa$. By the facts quoted above, let $\mu>\kappa$ be such that $\mathcal{P}(\kappa)/\mathcal{I}$ is isomorphic to the Boolean algebra ${\rm Borel}({}^\mu 2)/\mathcal{J}$, where $\mathcal{J}$ is the ideal of null sets in ${}^\mu 2$. We fix an isomorphism $\jmath$ which exemplifies this fact. Fix any cardinal $\theta<\kappa$.

Viewing the Boolean algebra as a forcing notion, let $\name{\eta}=(\name{\eta}_\alpha:\alpha<\mu)$ be a random sequence of reals (i.e. $\name{\eta}_\alpha\in{}^\omega 2$ for every $\alpha<\mu$, and $\name{\eta}_\alpha=\langle\name{\eta}(\omega\alpha+n):n\in\omega\rangle$). For each $\name{\eta}_\alpha$ we choose a sequence of sets $(B_{\alpha n}:n\in\omega)$ so that $B_{\alpha n}\subseteq\kappa$ and $\jmath(B_{\alpha n}/\mathcal{I})=({}^\mu 2)^{[(\alpha,\name{\eta}_\alpha(n))]}/ \mathcal{J}$, i.e. $\jmath(B_{\alpha n}/\mathcal{I})=\{\nu\in{}^\mu 2:\nu(\alpha)=\name{\eta}_\alpha(n)\}/ \mathcal{J}$.

For every $\alpha<\mu$ and each $n\in\omega$ we define $e_{\alpha n}\in{}^\kappa 2$ by $e_{\alpha n}(i)=1 \Leftrightarrow i\in B_{\alpha n}$. 
We concentrate on the collection $T=\{e_{\alpha n}:\alpha<\theta,n\in\omega\}$.
We define $\kappa$ colorings $c_i$ for every $i<\kappa$, each $c_i$ is a function from $\theta\times\omega$ into $2$, by letting $c_i(\alpha,n)=e_{\alpha n}(i)$ for each member of $T$. We claim that for some $i<\kappa$ the coloring $c_i$ exemplifies the negative relation $\binom{\theta}{\omega}\nrightarrow\binom{\aleph_1}{\omega}^{1,1}_2$.

For proving this, choose an ordinal $i<\kappa$ for which $\{(e_{\alpha n}(i):n\in\omega):\alpha<\theta\}$ is a Sierpi\'nski set, i.e. for every $B\subseteq{}^\omega 2$ of Lebesgue measure zero we have $(e_{\alpha n}(i):n\in\omega)\notin B$ apart from a countable set of such sequences. 
Notice that each sequence of the form $(e_{\alpha n}(i):n\in\omega)$ is an element in ${}^\omega 2$.
Such $i$ exists by the following paragraph, upon noticing that $\theta<\kappa$.

The existence of a Sierpi\'nski set is a well-known property of real-valued measurable cardinals, see \cite{MR1234282} 6F, p. 215.
We make the comment that here is the only point along the proof in which we use the assumption $\theta<\kappa$, and we do not know whether the existence of a Sierpi\'nski set can be guaranteed if $\theta=\kappa$.

Assume now that $H_0\in[\theta]^{\aleph_1}, H_1\in[\omega]^{\aleph_0}$, and we shall show that $c_i\upharpoonright(H_0\times H_1)$ is not constant. Let $B$ be $\{e\in{}^\omega 2: e\upharpoonright H_1$ is constant$\}$. Since $H_1$ is unbounded in $\omega$ we have $m(B)\leq\frac{1}{2^n}$ for every $n\in\omega$ and hence $m(B)=0$. Consequently, $(e_{\alpha n}(i):n\in\omega)\notin B$ apart from a countable set of such sequences.
Since $|H_0|>\aleph_0$ we can choose an ordinal $\alpha\in H_0$ so that $(e_{\alpha n}(i):n\in\omega)\upharpoonright H_1$ is not constant. It means that there are $n_0,n_1\in H_1$ such that $c_i(\alpha,n_0)\neq c_i(\alpha,n_1)$, so we are done.

\hfill \qedref{nnegbelowrvm}

As a consequence, we deduce the following known fact:

\begin{corollary}
\label{ssandrrvm} If there is a real-valued measurable cardinal, then the splitting number is $\aleph_1$.
\end{corollary}

\par\noindent\emph{Proof}. \newline 
By Claim \ref{belowssss}, if $\mathfrak{s}>\aleph_1$ then $\binom{\omega_1}{\omega}\rightarrow\binom{\omega_1}{\omega}^{1,1}_2$. However, this is impossible if one assumes the existence of a real-valued measurable cardinal. 

\hfill \qedref{ssandrrvm}

Finally, we return to Question \ref{rrvm} and we ask about the negative relation $\binom{\kappa}{\omega}\nrightarrow\binom{\kappa}{\omega}^{1,1}_2$ for the real-valued measurable $\kappa$ itself. If we force the existence of such a cardinal by the classical way of Solovay, \cite{MR0290961}, then $\mathfrak{r}=\kappa$ and the negative relation follows. It turns out that this is always true:

\begin{claim}
\label{mmc} If $\kappa$ is real-valued measurable then $\mathfrak{r}\geq\kappa$. Consequently, if $\mathfrak{c}=\kappa$ is real-valued measurable then $\binom{\kappa}{\omega}\nrightarrow\binom{\kappa}{\omega}^{1,1}_2$.
\end{claim}

\par\noindent\emph{Proof}. \newline 
Suppose that $\kappa$ is real-valued measurable and assume toward contradiction that $\mathfrak{r}<\kappa$. Since $\kappa$ is a limit cardinal, $\theta=\mathfrak{r}^+<\kappa$ as well. By Theorem \ref{nnegbelowrvm} we have $\binom{\theta}{\omega}\nrightarrow\binom{\theta}{\omega}^{1,1}_2$. However, $\binom{\theta}{\omega}\rightarrow\binom{\theta}{\omega}^{1,1}_2$ by Theorem \ref{aboverrrr}, a contradiction.

In case $\mathfrak{c}=\kappa$ is real-valued measurable we have $\mathfrak{r}=\kappa$, hence $\binom{\kappa}{\omega}\nrightarrow\binom{\kappa}{\omega}^{1,1}_2$ by virtue of Corollary \ref{ccc}, so we are done.

\hfill \qedref{mmc}

\newpage

\bibliographystyle{amsplain}
\bibliography{arlist}

\end{document}